\documentclass[10pt,twoside]{article}
\usepackage{amssymb}
\usepackage{amsmath}
\usepackage{latexsym}
\newcommand{\COM}[1]{}
\renewcommand{\theequation}{\arabic{section}.\arabic{equation}}

\newtheorem{theorem}{Theorem}[section]

\newtheorem{lemma}{Lemma}[section]

\newtheorem{corollary}{Corollary}[section]

\newtheorem{remark}{\normalfont\scshape Remark}[section]

\newenvironment{proof}{\noindent\textsc{Proof.\/}}{}

\newcommand{\subj}[2]{\textsf{AMS 2000 subject classifications.}
Primary {#1}; Secondary {#2}.\newline}
\newcommand{\key}[1]{\textsf{Keywords and phrases.} {#1}.\newline}
\newcommand{\abb}[1]{\textsf{Abbreviated title.} {#1}.}

\newcommand{\fot}[5]{\renewcommand\thefootnote{}
\footnotetext{\parindent=0.0mm \vskip-3mm \subj{#1}{#2}\key{#3}\abb{#4}
\newline\textsf{Date.} \date{\today}}}

\COM{%%%%%%%%%%%%%%%%%%%%%%%%%%%%%%

}%%%%%%%%%%%%%%%%%%%%%%%%%%%%%%%%%%%%%%%

\def\vsb{\hfill$\Box$}

\def\vsp{\vskip-8mm\hfill$\Box$\vskip3mm}

\newcommand{\be}{\begin{equation}}
\newcommand{\ee}{\end{equation}}
\newcommand{\bea}{\begin{eqnarray}}
\newcommand{\eea}{\end{eqnarray}}
\newcommand{\beaa}{\begin{eqnarray*}}
\newcommand{\eeaa}{\end{eqnarray*}}
\newcommand{\beal}{\begin{aligned}}
\newcommand{\eeal}{\end{aligned}}

\newcommand{\var}{\mathrm{Var\,}}

\newcommand{\ifff}{\ttt{$\Longleftrightarrow$}}

\newcommand{\bn}{\mathbf n}

\newcommand{\bk}{\mathbf k}

\newcommand{\sumk}{\sum^n_{k=1}}

\newcommand{\ttt}[1]{\quad\mbox{ #1}\quad}

\newcommand{\asto}{\stackrel{a.s.}{\to}}

\newcommand{\nifi}{n\to\infty}

\newcommand{\iid}{i.i.d.\ }

\newcommand{\zpd}{Z\hskip-7pt Z_{+}^d}

\newcommand{\window}{\sum_{j=n+1}^{n+k}X_j}
\newcommand{\tna}{T_{n,n+n^\alpha}}

\newcommand{\tnan}{T_{n,n+a_n}}

\newcommand{\lsv}{L\in\cal{SV}}

\newcommand{\fenster}{T_{\bn,\bn+\bna}}

\newcommand{\ba}{\boldsymbol{\alpha}}
\newcommand{\bna}{\bn^{\ba}}
\newcommand{\bbna}{|\bna|}

\newcommand{\bbn}{|\bn|}

\begin{document}
%\tableofcontents
\date{}
\title{\textsf{On the LSL for random fields}}
\author{Allan Gut\\Uppsala University \and  Ulrich Stadtm\"uller\\
University of Ulm}
\maketitle

\begin{abstract}\noindent
In some earlier work we have considered extensions of Lai's (1974) law
of the single logarithm for delayed sums to a multiindex setting with
the same as well as different expansion rates in the various
dimensions. A further generalization concerns window sizes that are
regularly varying with index 1 (on the line).  In the present paper we
establish multiindex versions of the latter as well as for some mixtures
of expansion rates. In order to keep things within reasonable size we
confine ourselves to some special cases for the index set $\mathbb{Z}_+^2$.
\end{abstract}

\fot{60F15, 60G50}{60G70}{Delayed sums, window, law of the iterated
logarithm, law of the single logarithm,
sums of \iid random variables, window, slowly varying function,
multidimensional indices, random fields}{On the LSL for random fields}

\markboth{A.\ Gut and U.\ Stadtm\"uller}{The LSL for random fields}

\section{Introduction}
\setcounter{equation}{0}

Let $X,\,\{X_k,\, k\geq1\}$ be \iid random variables with mean
$\mu=0$ and partial sums $\{S_n,\,n\geq1\}$.  The Hartman-Wintner
\emph{Law of the Iterated Logarithm\/} (LIL) \cite{hw} states that
\[
\limsup_{\nifi}
\frac{S_n}{\sqrt{2 n \log \log n}}
=\sigma \mbox{\quad a.s.}\]
Later Strassen, in \cite{strassen}, proved the necessity in the sense
that he showed that if
\[P\Big(\limsup_{\nifi}
\frac{|S_n|}{\sqrt{n \log \log n}}<\infty \Big)>0\,,\]
then $E\,X^2 <\infty$ and $E\,X=0$.
\begin{remark} \emph{By the Kolmogorov zero-one law it follows that if
the probability of the limsup being finite is positive, then it is
automatically equal to 1.}\end{remark}
\begin{remark} \emph{Here and throughout there also exist obvious analogs
for the limit inferior}.\vsb
\end{remark}
The \emph{Law of the Single Logarithm\/} (LSL) is due to Lai
\cite{lai}, and deals with \emph{delayed sums\/} or \emph{windows},
viz.,
\[T_{n,n+k}=\window,\quad n\geq0,\, k\geq 1,\]
and states that, if
\[E\,|X|^{2/\alpha}\big(\log^+|X|\big)^{-1/\alpha}<\infty\ttt{and}
E\,X=0,\;E\,X^2=\sigma^2,
\]
then, for  $0<\alpha<1$,
\[ \limsup_{\nifi}\frac{\tna}{\sqrt{2n^\alpha\log n}}=\sigma\,\sqrt{1-\alpha}
\mbox{\quad a.s.},
\]
and, conversely, that if
\[P\Big(\limsup_{\nifi}
\frac{|\tna|}{\sqrt{n^\alpha\log n}}<\infty \Big)>0\,,\]
then
\[E\,|X|^{2/\alpha}\big(\log^+|X|\big)^{-1/\alpha}<\infty\ttt{and}E\,X=0.
 \]
\begin{remark}\emph{Here, and throughout, $\log^+ x=\max\{\log x,1\}$
for $x>0\,.$}\vsb
\end{remark}
These results can be generalized in various ways. However, let us
first mention that there also exist \emph{one-sided\/} versions of the
above results.

Martikainen \cite{mart}, Rosalsky \cite{ros} and Pruitt
\cite{pruitt} independently proved a one-sided LIL to the effect that
if
\[-\infty<\limsup_{\nifi}
\frac{S_n}{\sqrt{2 n \log \log n}}<\infty\]
with positive probability, then $E\,X^2< \infty$ and $E\,X=0$. An
analogous one-sided version of Lai's result is Theorem 3 of
\cite{cl75}, where now the finiteness of the one-sided limsup is
equivalent to the same moment condition as in \cite{lai}, however,
based on $X^+=\max\{X,0\}$ rather than on $|X|$.

As for multiindex results,
Wichura \cite{wichura} proved the following LIL for random fields:
Let $\{X_{\bk},\,\bk\in \zpd\}$ be \iid random variables with
partial sums $S_{\bn} =\sum_{\bk\leq \bn}X_{\bk}$, $\bn\in\zpd$, where
the random field or index set $\zpd$, $d\geq 2$, is the positive
integer $d$-dimensional lattice with coordinate-wise partial ordering
$\leq $, where $\bbn=\prod_{k=1}^{d}n_{k}$, and where $\mathbf{n\to
\infty}$ means that $ n_{k}\to \infty $ for all $k=1,2,\ldots
,d$. Then
\[\limsup_{\bn \to \infty} \, \frac{S_\bn}{\sqrt{2\bbn
\log \log\bbn}} =\sigma\sqrt{d}\ttt{a.s.,} \]
provided
\[E\,X^2\frac{(\log^+|X|)^{d-1}}{\log^+\log^+|X|}
<\infty \ttt{and} E\,X=0,\;  E\,X^2=\sigma^2\,,\]
and, conversely, if
\[P\Big(\limsup_{\bn \to \infty} \, \frac{|S_\bn|}{\sqrt{2\bbn
\log \log\bbn}}<\infty\Big)>0\,,\]
then
\[E\,X^2\frac{(\log^+|X|)^{d-1}}{\log^+\log^+|X|}
<\infty \ttt{and} E\,X=0.\]
The analogous multiindex extension of Lai's result is given in
\cite{gs1}. In order to describe the main result there, let
$\{X_{\bk},\,\bk\in \zpd\}$ be \iid random variables and define
the delayed sums---which in this setting turn into ``real'' windows,
\[T_{\bn,\bn+\bn^{\alpha}}= \sum_{\bn\le \bk \le \bn +\bn^\alpha} X_\bk ,
\, \bn\in\zpd\,,
\]
and where addition is to be taken coordinate-wise. Then, for
$0<\alpha<1$,
\[\limsup_{\bn\to\infty}\,
\frac{T_{\bn,\bn+\bn^{\alpha}}}{\sqrt{2|\bn|^{\alpha}\log\bbn}}
=\sigma\sqrt{1-\alpha} \ttt{a.s.,}\]
provided
\[E\,|X|^{2/\alpha}(\log^+|X|)^{d-1-1/\alpha}<\infty\ttt{and}E\,X=0,\;
E\,X^2=\sigma^2.\]
We also proved that if
\[P\Big(\limsup_{\bn \to \infty} \,
\frac{|T_{\bn,\bn+\bn^{\alpha}}|}{\sqrt{\bn|^{\alpha}\log\bbn}}<\infty
\Big)>0\,,\]
then
\[E\,|X|^{2/\alpha}(\log^+|X|)^{d-1-1/\alpha}<\infty\ttt{and}E\,X=0.\]
Note in particular that the moment condition depends on
the $d$ (as in Wichura's LIL above).

In \cite{gs2} this result was generalized to allow for different $\alpha$'s
for the different directions, and it was shown there that if
\[E\,|X|^{2/\alpha_1}(\log^+|X|)^{p-1-1/\alpha_1}<\infty
\ttt{and}E\,X=0,\; E\,X^2=\sigma^2,\]
then
\[\limsup_{\bn\to\infty}\,
\frac{\fenster}{\sqrt{2\bbna\log\bbn}}
=\sigma\sqrt{1-\alpha_1} \ttt{a.s.,}
\]
and, conversely, that if
\[P\Big(\limsup_{\bn\to\infty}\,
\frac{|\fenster|}{\sqrt{\bbna\log\bbn}}<\infty\Big)>0\,,\]
then
\[E\,|X|^{2/\alpha_1}(\log^+|X|)^{p-1-1/\alpha_1}<\infty
\ttt{and}E\,X=0,
\]
where now
\[
0<\alpha_1\leq\alpha_2\leq\cdots\leq\alpha_d<1 \,, \
p=\max\{k:\alpha_k=\alpha_1\}\ttt{and} \ \bbna=\Pi_{k=1}^d n_k^{\alpha_k}\,.
\]
\begin{remark}\emph{When $\alpha_1=\alpha_2=\cdots=\alpha_d=\alpha$ the
theorem reduces (of course) to the previous one}.\vsb
\end{remark}

The next natural question is to consider the boundary cases
$\alpha=0$ and $\alpha=1$. As for the former, one case is the trivial one in
which the windows reduce to single random variables, viz.,
$T_{n,n+1}=X_{n+1}$. A nondegenerate variation concerns the delayed sums
$T_{n,n+\log n}$, $n\geq1$, which obey the
so-called Erd\H os-R\'enyi law (\cite{er70}, Theorem 2, \cite{cr81},
Theorem 2.4.3).

In \cite{gjs} we considered the nondegenerate boundary case  $\alpha=1$
at the other end. In this case the  window size is larger
than any power less than one, and at the same time not quite
linear. Technically, the paper focused on windows of the form
\bea\label{an} \tnan\ttt{where} a_n&=&\frac{n}{L(n)} \ttt{with}\\
\label{an1} \ttt{a differentiable function } L(\cdot)\nearrow \infty
\in \cal{SV} &\mbox{and}& \frac{xL'(x)}{L(x)}\searrow \ttt{as}
x\to\infty\,, \eea where $\lsv$ means that $L$ is slowly varying at
infinity (see e.g.\ \cite{bgt} or \cite{g07}, Section A.7), and where,
for convenience, we shall permit ourselves to treat quantities such as
$a_n=n/L(n)$ and $a_n=n/\log n$ and so on as integers.

Finally, let
\[d_n=\log\frac{n}{a_n} +\log \log n=\log L(n)+\log \log n,\quad n\geq2,\]
and set
\[f_n = \min\{a_n\cdot d_n,\,n\},\]
with $f$ an increasing interpolating function, i.e., $f(x)=f_{[x]}$
for $x>0$, and with $f^{-1}$ the corresponding (suitably defined)
inverse function. Then, in short, (the precise equivalence has to be
formulated as above)
\[\limsup_{\nifi}\,\frac{\tnan}{\sqrt{2a_n d_n}}
=\sigma \ttt{a.s.}
\ifff E\big(f^{-1}(X^{2})\big)<\infty,\;E\,X=0,\;
E\,X^2=\sigma^2.
\]
As an introduction and point of departure to what follows, here are
the two canonical examples; \cite{gjs}, Corollaries 2.1 and 2.2, which
concern the cases $L(n)=\log n$ and $L(n)=\log\log n$, that is,
\[
d_n=2\log\log n,\qquad f_n=2\frac{n}{\log n}\log\log n\,,
\]
and
\[d_n=\log\log\log n+\log\log n\sim\log\log n,\qquad f_n\sim n\,,
\]
respectively.

\begin{theorem}\label{thm75} Let $X,\,\{X_k,\, k\geq1\}$ be \iid random
variables with partial sums $S_n=\sumk X_k$, $n\geq1$.\\\noindent
\emph{(i)} If
\[
E\,X^2\frac{\log^+|X|}{\log^+\log^+|X|}<\infty
\ttt{and}E\,X=0,\;\var X=\sigma^2,
\]
then
\[\limsup_{\nifi}\,\frac{T_{n,n+n/\log n}}
{\sqrt{4\frac{n}{\log n}\log\log n}}=\sigma\ttt{a.s.}
\]
Conversely, if
\[P\Big(\limsup_{\nifi}\,\frac{|T_{n,n+n/\log n}|}
{\sqrt{\frac{n}{\log n}\log\log n}}<\infty\Big)>0\,,\]
then $E\,X^2\frac{\log^+|X|}{\log^+\log^+|X|}<\infty$ and $E\,X=0$.\\[1.4mm]
\noindent
\emph{(ii)} If  $E\,X=0$ and $E\,X^2=\sigma^2<\infty$, then
\[\limsup_{\nifi}\,\frac{T_{n,n+n/\log\log n}}
{\sqrt{2n}}=\sigma\ttt{a.s.}
\]
Conversely, if
\[P\Big(\limsup_{\nifi}\,\frac{|T_{n,n+n/\log\log n}|}
{\sqrt{n}}<\infty\Big)>0\,,\]
then   $E\,X^2<\infty$ and $E\,X=0$.  \vsb
\end{theorem}

The purpose of the present paper is to investigate random field
analogs of this result. However, in order to keep things within
reasonable bounds we shall confine ourselves to the case $d=2$ and
the windows
\begin{equation}\label{window}
T_{(m,n)\,,\,(m+a^{(1)}_m,n+a^{(2)}_n)}
=\sum_{\stackrel{m \le i \le m+a^{(1)}_m}{n \le j \le n+a^{(2)}_n}}
X_{i,j}\,,\end{equation}

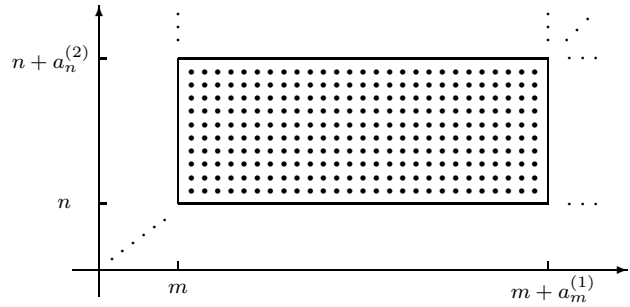
\begin{figure} [ht]
\begin{center}
\begin{picture}(200,90)(0,0)
\put(0,0){\vector(1,0){210}}\put(10,-10){\vector(0,1){110}}

\multiput(15,4)(4,3){6}{\circle*{1}}\multiput(187,87)(4,4){3}{\circle*{1}}
\multiput(188,25)(5,0){3}{\circle*{1}}\multiput(188,80)(5,0){3}{\circle*{1}}
\multiput(40,87)(0,5){3}{\circle*{1}}\multiput(180,87)(0,5){3}{\circle*{1}}
\put(40,25){\line(1,0){140}}\put(40,80){\line(1,0){140}}
\put(40,25){\line(0,1){55}}\put(180,25){\line(0,1){55}}

\put(40,0){\line(0,1){3}}\put(180,0){\line(0,1){3}}

\put (40,-7){\makebox(0,0){\footnotesize $m$}}
\put (182,-7){\makebox(0,0){\footnotesize $m+a^{(1)}_m$}}

\put(10,25){\line(1,0){3}}\put(10,80){\line(1,0){3}}
\put (-3,25){\makebox(0,0){\footnotesize $n$}}
\put (-8,80){\makebox(0,0){\footnotesize $n+a^{(2)}_n$}}
\multiput(45,30)(5,0){27}{\circle*{2}}
\multiput(45,40)(5,0){27}{\circle*{2}}
\multiput(45,50)(5,0){27}{\circle*{2}}
\multiput(45,60)(5,0){27}{\circle*{2}}
\multiput(45,70)(5,0){27}{\circle*{2}}
\multiput(45,35)(5,0){27}{\circle*{2}}
\multiput(45,45)(5,0){27}{\circle*{2}}
\multiput(45,55)(5,0){27}{\circle*{2}}
\multiput(45,65)(5,0){27}{\circle*{2}}
\multiput(45,75)(5,0){27}{\circle*{2}}

\end{picture}\end{center}
\caption{\label{belieb} A typical window.}
\end{figure}
\noindent where
\begin{enumerate}
\item $a^{(1)}_m=m/\log m$ and  $a^{(2)}_n=n/\log n$,
\item $a^{(1)}_m=m/\log\log m$ and  $a^{(2)}_n=n/\log\log n$,
\item $a^{(1)}_m=m/\log m$ and  $a^{(2)}_n=n/\log\log n$,
\end{enumerate}
leaving further cases to the reader(s).

In Section \ref{setup} we present our setup and main result, after
which some preliminaries are given in Section \ref{prel}. Following
this, we provide proofs in Section \ref{pfthm}.
In Section \ref{7075} we connect our
earlier papers in the area in that we consider windows of the form
\[T_{(m,n)\,,\,(m+m^{\alpha},\,n+n/\log n)}.\]
Some final remarks, preceding an Appendix in which we  have collected
some technicalities, are given in Section \ref{anm}.

\section{Setup and main result}\label{setup}
\label{main}\setcounter{equation}{0} In the remainder of the paper
we thus suppose that $X$ and $\{X_{i,j},\, i,j\geq1\}$ are \iid
random variables with partial sums
$S_{m,n}=\sum_{(i,j)=(1,1)}^{(m,n)}X_{i,j}$ and windows
$T_{(m,n),\,(m+a^{(1)}_m,n+a^{(2)}_n)}$ (see (\ref{window})) for
$m,n\geq1$.
\begin{theorem} \label{thmmain} \emph{(i)} If
\[E\,X^2\frac{(\log^+|X|)^3}{\log^+\log^+|X|}<\infty
\ttt{and}E\,X=0,\;E\,X^2=\sigma^2\,,
\]then
\[
\limsup_{m,n\to\infty}\,\frac{T_{(m,n)\,,\,(m+m/\log m,n+n/\log n)}}
{\sqrt{4mn\frac{\log\log m + \log\log n}{\log m\log n}}}=\sigma\ttt{a.s.}
\]
Conversely, if
\[P\Big(\limsup_{\nifi}\,\frac{|T_{(m,n)\,,\,(m+m/\log m,n+n/\log n)}|}
{\sqrt{mn\frac{\log\log m + \log\log n}{\log m\log n}}}<\infty\Big)>0\,,\]
then $E\,X^2\frac{(\log^+|X|)^3}{\log^+\log^+|X|}<\infty$ and $E\,X=0$.\\[1mm]
\noindent
\emph{(ii)} If
\[E\,X^2\log^+|X|\,\log^+\log^+|X|<\infty\ttt{and}E\,X=0,\;E\,X^2=\sigma^2\,,
\]
then
\[
\limsup_{m,n\to\infty}\, \frac{T_{(m,n)\,,\,(m+m/\log\log m,n+n/\log\log n)}}
{\sqrt{2mn\frac{\log\log m +\log\log n}{\log\log m \log\log n}}}
=\sigma\ttt{a.s.}\]
Conversely, if
\[
P\Big(\limsup_{\nifi}\,\frac{|T_{(m,n)\,,\,(m+m/\log\log m,n+n/\log\log n)}|}
{\sqrt{mn\frac{\log\log m+\log\log n}{\log\log m \log\log n}}}<\infty\Big)>0\,,
\]
then $E\,X^2\log^+|X|\,\log^+\log^+|X|<\infty$ and $E\,X=0$.\\[1mm]
\noindent
\emph{(iii)} If
\[E\,X^2 (\log^+|X|)^2<\infty\ttt{and}E\,X=0,\;E\,X^2=\sigma^2\,,\]
then
\[\limsup_{m,n\to\infty}\, \frac{T_{(m,n)\,,\,(m+m/\log m,n+n/\log\log n)}}
{\sqrt{4mn\frac{\log\log m+\log\log n}{\log m\log\log n}}}
=\sigma\ttt{a.s.}\]
Conversely, if
\[P\Big(\limsup_{\nifi}\,\frac{|T_{(m,n)\,,\,(m+m/\log m,n+n/\log\log n)}|}
{\sqrt{4mn\frac{\log\log m+\log\log n}{\log m\log\log n}}}
<\infty\Big)>0\,,\]
then $E\,X^2 (\log^+|X|)^2<\infty$ and $E\,X=0$.
\end{theorem}

\section{Preliminaries}
\label{prel}\setcounter{equation}{0}
Before we jump into the proof of Theorem \ref{thmmain} we present the
Kolmogorov exponential bounds (cf.\ e.g.\ \cite{g07}, Section 8.2)
adapted to the present situation, a lemma that relates certain tail
sums to moments together with some special cases that will be of use
later.

But first some notation.
Let, for $i=1,2$ and $n\geq2$,
\[d_n^{(i)}=\log\frac{n}{a_n^{(i)}}+\log \log n=\log L_i(n)+\log \log n
\ttt{and}r_n=\max\{d_n^{(1)},\, d_n^{(2)}\}.\]
The two-dimensional analogs of $a_n$, $d_n$, and $f_n$ are
as follows:
\[a_{m,n}=a_m^{(1)}\cdot a_m^{(2)}=\frac{m}{L_1(m)}\cdot\frac{n}{L_2(n)},\qquad
d_{m,n}=r_m + r_n,\qquad f_{m,n}=a_{m,n}\cdot d_{m,n}\,,\] where in
the following the slowly varying functions $L$ are
logarithms or iterated logarithms.

In order to introduce the traditional double-truncation we let $\delta>0$
and $\varepsilon>0$, and set
\bea\label{bn} b_{m,n}
 =\frac{\sigma\delta}{\varepsilon} \sqrt{\frac{a_{m,n}}{d_{m,n}}},
 \eea
and note that \bea\label{adrotf}
\frac{a_{m,n}}{b_{m,n}}=\frac{\varepsilon}{\sigma\delta}\sqrt{f_{m,n}}\,.
\eea Next we set \beaa X_{m,n}'&=&X_{m,n}I\{|X_{m,n}|\leq
b_{m,n}\},\quad
X_{m,n}''=X_{m,n}I\{b_{m,n}<|X_{m,n}|< \delta\sqrt{f_{m,n}}\}, \\[3pt]
X_{m,n}'''&=&X_{m,n}I\{|X_{m,n}|\geq\delta\sqrt{f_{m,n}}\}\,. \eeaa
Throughout the following all objects with primes or multiple primes
refer to the respective truncated summands.

\subsection{Exponential bounds}
\label{kolmexp}

The standard procedure for estimating $E\,X_{i,j}'$ yields
\[|E\,X_{i,j}'|\leq E|X|I\{|X|> b_{i,j}\}\leq\frac{E\,X^2I\{|X|>b_{i,j}\}}
{b_{i,j}},\]
so that, omitting intermediate steps and recalling (\ref{adrotf}),
\beaa
|E\,T_{(m,n)\,,\,(m+a_m^{(1)},n+a_n^{(2)})}'|
&\leq& a_{m,n}\frac{E\,X^2I\{|X|>b_{m,n}\}}{b_{m,n}}\\
&=&\frac{\varepsilon}{\sigma\delta} \sqrt{f_{m,n}}
  E\,X^2I\{|X|>b_{m,n}\}\\[3pt]
&=&o(\sqrt{f_{m,n}})\ttt{as}m,n\to\infty. \eeaa
Moreover, a simple calculation yields, for any $\delta>0$ and
$\max\{m,n\}\ge n_0(\delta)$,
\[a_{m,n}\sigma^2(1-\delta)\leq
\var \big(T_{(m,n)\,,\,(m+a_m^{(1)},n+a_n^{(2)})}' \big)\leq a_{m,n}\sigma^2.\]
Inserting the above estimates into the Kolmogorov exponential bounds then
yields
\bea
P\Big(T_{(m,n)\,,\,(m+a_m^{(1)},n+a_n^{(2)})}'>\varepsilon\sqrt{2f_{m,n}}\Big)
\;\begin{cases}&
\leq \exp\big\{-\frac{\varepsilon^2(1-\delta)^3}{\sigma^2}
\cdot d_{m,n}\big\}\,,\\[6pt]
&\geq \exp\big\{-\frac{\varepsilon^2(1+\delta)^2(1+\gamma)}{\sigma^2(1-\delta)}
\cdot d_{m,n} \big\}\,,
\end{cases}\label{expbds}
\eea for any $\gamma>0$ and $\max\{m,n\}$ large.

\subsection{A useful lemma}
\label{lemma} The present subsection contains a technical lemma and some
consequences that will be of use later.
\begin{lemma}\label{lemmam} Consider a positive function G on
$[1,\,\infty)\times [1,\,\infty]$,
such that, for some constants $0<c_1\le c_2$ and for all
$m,n \in \mathbf{N}$,
\begin{equation} \label{cc} c_1\le \frac{G(m,n)}{G(u,v)}\le c_2 \ttt{for}
(u,v)\in [m,\, m+1]\times[n,\,n+1]\,.
\end{equation}
Furthermore, let
\begin{equation}
M(x)\,=\,\ \big|\{(u,v)\in [1,\, \infty)\times [1,\,\infty):
G(u,v)\le x\}\big|\,,
\end{equation}
and assume that, for any $c>1$, the quotient $\frac{M(cx)}{M(x)}$ is
bounded as $x \to \infty\,.$ Finally, let $X$ be a nonnegative
random variable. Then
\[\sum_{m,n=1}^\infty P\big(X> G(m,n)\big)<\infty
\ifff \iint P\big(X>G(u,v)\big)\,dudv<\infty \ifff
E\,M(X)<\infty\,.\]
\end{lemma}
\begin{proof} The equivalences follow from
\[ \sum_{m,n=1}^\infty P\big(X> G(m,n)\big)
= \sum_{m,n=1}^\infty \int_m^{m+1}\int_n^{n+1}
P\Big(X> \frac{G(m,n)}{G(u,v)}G(u,v)\Big)du dv\,,\]
(\ref{cc}), and the fact that
\beaa \iint P\big(X>G(u,v)\big)\,dudv
&=&\iint\Big(\int_{G(u,v)}^{\infty} dF_X(x)\Big)du\,dv\\
&=&\int_0^\infty\Big(\iint_{\{(u,v):G(u,v)\leq x\}}\,dudv\Big) dF_X(x)\\[3pt]
&=&\int_0^\infty M(x)\, dF_X(x)=E\, M(X). \eeaa\vsp
\end{proof}

Next, some useful examples that we collect separately. The
first part  is immediate from the lemma and the Borel--Cantelli
lemmas. For the verification of the special cases we refer to the
appendix.
\begin{corollary}\label{corrmom} Suppose that $\{X\,,\;X_{m,n},\,m,n\geq1\}$
are \iid random variables and let $M$ and $G$ be defined as in the lemma. Then
\[\limsup_{m,n\to\infty}\frac{|X_{m,n}|}{\sqrt{G(m,n)}}<\infty\mbox{\quad a.s.}
\ifff\sum_{m,n}P(|X|>\sqrt{G(m,n)})<\infty\ifff
E\,M(X^2)<\infty\,.\]
In particular,
\bea
\limsup_{m,n\to\infty}\frac{|X_{m,n}|}{\sqrt{mn\frac{\log\log m+\log \log n}
{\log m\log n}}}<\infty\quad\mbox{a.s.}
&\iff&\sum_{m,n}P\Big(|X|>\sqrt{mn\frac{\log\log m+\log \log n}
{\log m\log n}}\Big)<\infty\nonumber\\
&\iff& E\,X^2\frac{(\log^+|X|)^3}{\log^+\log^+|X|}<\infty\,,
\label{corrlog}\\[4pt]
\limsup_{m,n\to\infty}\frac{|X_{m,n}|}{\sqrt{mn\frac{\log\log
m+\log\log n}{\log\log m \,\log\log n}}}<\infty\quad\mbox{a.s.}
&\iff&\sum_{m,n}P\Big(|X|>\sqrt{mn\frac{\log\log m+\log\log n}
{\log\log m \,\log\log n}}\Big)<\infty\nonumber\\
&\iff&  E\,X^2\log^+|X|\log^+\log^+|X|<\infty\,,\label{corrloglog}\\[4pt]
\limsup_{m,n\to\infty}\frac{|X_{m,n}|}{\sqrt{mn\frac{\log\log m
+\log \log n}{\log m \, \log \log n}}}<\infty\quad\mbox{a.s.}
&\iff&\sum_{m,n}P\Big(|X|>\sqrt{mn\frac{\log\log m+\log \log n}{\log
m \log \log n}}\;\Big)
<\infty\nonumber\\
&\iff& E\,X^2 (\log^+|X|)^2<\infty\,. \label{corrlogloglog} \eea
\end{corollary}

\section{Proof of Theorem \ref{thmmain}}
\setcounter{equation}{0}\label{pfthm}
Since the proof follows the pattern of our predecessors
we confine ourselves to providing the proof of
part (i) of the theorem in somewhat more detail and the other cases more
sketchily.

\subsection{Proof of Theorem \ref{thmmain}(i)}

\subsubsection*{Exponential bounds}
In this case we have
\[a_{m,n}=\frac{mn}{\log m\log n},\qquad
d_{m,n}=2(\log\log m+\log\log n),\qquad f_{m,n}=\frac{2mn}{\log
m\log n}\cdot(\log\log m+\log\log n)\,,\] so that the  Kolmogorov
exponential bounds (\ref{expbds}) yield \bea &&\hskip-5pc
P\Big(T_{(m,n)\,,\,(m+m/\log m,n+n/\log n)}'
>\varepsilon\sqrt{2\frac{2mn(\log\log m+\log\log n)}{\log m\log n}}\;\Big)
\nonumber\\
&&\hskip3pc\begin{cases}&
\leq \exp\big\{-\frac{\varepsilon^2(1-\delta)^3}{\sigma^2}
\cdot 2(\log\log m+\log\log n)\big\}\,,\\[6pt]
&\geq
\exp\big\{-\frac{\varepsilon^2(1+\delta)^2(1+\gamma)}{\sigma^2(1-\delta)}
\cdot 2(\log\log m+\log\log n)\big\}\,,
\end{cases}\label{expbdslog}
\eea for any $\gamma>0$ and $\max\{m,n\}$ large.

\subsubsection*{Upper bound}

In order to obtain an upper bound we consider, for some constant
$c>0$ to be chosen later, the subset
\bea\label{minj}
(m_i,n_j)=\big(e^{\sqrt{c \,i}},e^{\sqrt{c\,j}}\big)\ttt{for}i,j
\geq \max\{1/c\,,\,1\}\,.
\eea
The upper exponential bound then reduces to
\bea
&&\hskip-5pc P\Big(T_{(m_i,n_j)\,,\,(m_i+m_i/\log
m_i,n_j+n_j/\log n_j)}'>\varepsilon
\sqrt{2\cdot\frac{2mn(\log\log m_i+\log\log n_j)}{\log m_i\log n_j}}\;\Big)
\nonumber\\
&&\hskip2pc\leq
C\,\exp\big\{-\frac{\varepsilon^2(1-\delta)^3}{\sigma^2} \cdot(\log
i+\log j)\big\} =C
(ij)^{\,-\,\frac{\varepsilon^2(1-\delta)^3}{\sigma^2}}\,,\label{expupplog}
\eea from which it follows that \bea
\sum_{i,j}P\Big(T_{(m_i,n_j)\,,\,(m_i+m_i/\log m_i,\,n_j+n_j/\log n_j)}'
>\varepsilon\sqrt{2\frac{2m_in_j(\log\log m_i+\log\log n_j)}
{\log m_i\log n_j}}\;\Big)<\infty \label{teins} \eea
whenever $\varepsilon>\sigma(1-\delta)^{-3/2}$.

As for the contribution of $T''$ we observe as in \cite{gs1} that
in order for $|T_{(m,n)\,,\,(m+a_{m}^{(1)},n+a_{n}^{(2)})}''|$
to surpass the level
$\eta\sqrt{f_{m,n}}$ it is necessary that at least
$N\geq\eta/\delta$ of the $X''$:s are nonzero, which, by stretching the
truncation bounds to the extremes, implies that
\bea
P(|T_{(m,n)\,,\,(m+a_{m}^{(1)},n+a_{n}^{(2)})}''|
>\eta\sqrt{f_{m,n}})
&\leq&\binom{a_{m}^{(1)}a_{n}^{(2)}}{N}\Big(P\big( b_{m,n}<|X|\leq\delta
\sqrt{f_{m,n}}\big)\Big)^{N}\nonumber\\
&\leq&C\Big(a_{m}^{(1)} a_{n}^{(2)}  \Big)^N
\bigg(\frac{E\,H(|X|)}{H( b_{m,n})}\bigg)^N\,,\label{ttzwei}
\eea
where $E\,H(|X|)<\infty$ is the appropriate moment condition.

In the present case this amounts, after simplifying, to
\bea
P''(m,n)&=&P(|T_{(m,n)\,,\,(m+m/\log m,n+n/\log n)}''|>\eta\sqrt{f_{m,n}})
\nonumber\\
&\leq& C\bigg(\frac{(\log\log m+\log\log n)\,(\log\log
(mn))}{(\log(mn))^3}\bigg)^N\,.\label{pzweiprel} \eea
This means that for our subset (\ref{minj}) we have
\bea\label{pzwei} P''(m_i,n_j)\leq C\bigg(\frac{(\log (ci)
+\log (cj))(\log(\sqrt{ci}+\sqrt{cj}))} {(\sqrt{ci}+\sqrt{cj})^3}\bigg)^N
\leq C\bigg(\frac{(\log i+\log j)(\log(i+j))}{i^{3/2}+j^{3/2}}\bigg)^N \,.
\eea
Since the sum of these
probabilities converges whenever $N\geq2$, we conclude, via the first
Borel--Cantelli lemma, considering in
addition that $N\delta\geq\eta$,  that
\bea\label{tzwei} \limsup_{i,j\to\infty}
\frac{|T_{(m_i,n_j)\,,\,(m_i+m_i/\log m_i,n_j+n_j/\log n_j)}''|}
{\sqrt{f_{m_i,n_j}}}\leq \delta \ttt{a.s.} \eea

The next step is to show that
\bea\label{tdrei}
\lim_{m,n\to\infty}\frac{T_{(m,n)\,,\,(m+m/\log m,n+n/\log n)}'''}
{\sqrt{f_{m,n}}}=0 \ttt{a.s.}
\eea
Now, since in order for $|T'''_{(m,n)\,,\,(m+m/\log m,n+n/\log n)}|$
to surpass the level $\eta\sqrt{f_{m,n}}$ infinitely often it is
necessary that infinitely many of the $X'''$:s are nonzero. However, via
an appeal to the first Borel--Cantelli lemma, the latter event has zero
probability. Namely, for every $\eta>0$ we have
\[
\sum_{m,n}P(|X_{m,n}|>\eta\sqrt{f_{m,n}}) =\sum_{m,n}
P\Big(|X|>\eta\sqrt{\frac{2mn}{\log m\log n}\cdot(\log\log
m+\log\log n)}\Big) <\infty\,,
\]
since, by assumption, $E\,X^2\frac{(\log^+|X|)^3}{\log^+\log^+|X|}<\infty$;
recall (\ref{corrlog}).

By combining (\ref{teins}), (\ref{tzwei}), and (\ref{tdrei}) we
are now in the position to conclude that \bea\label{zlog}
\limsup_{i,j\to\infty} \frac{T_{(m_i,n_j)\,,\,(m_i+m_i/\log
m_i,n_j+n_j/\log n_j)}} {\sqrt{2f_{m_i,n_j}}} \leq
\sigma(1-\delta)^{-3/2}+\delta \ttt{a.s.,} \eea for
$(m_i,n_j)=\big(e^{\sqrt{ci}},e^{\sqrt{cj}}\big)$ with $i,j \ge
1/c$ and any $c>0$.

\subsubsection*{Sufficiency for the entire field}

It remains to show that our process behaves accordingly for the
entire field. Assume for the moment that the random variables
are symmetric, let $\eta>0$, small, be given, and choose $c=\eta^2$.
Recalling that $a_n^{(1)}=a_n^{(2)}=n/\log n$, and that
$f_{m,n}=  2mn \frac{\log\log m +\log\log n}{\log m\log n}$, we thus consider
\begin{eqnarray*}\label{fillgap}\lefteqn{P\Big(\max_{\stackrel{m_i\le m\le
m_{i+1}}{n_j \leq n \leq n_{j+1}}}
\frac{T_{(m,n),(m+a_m^{(1)},n+a_n^{(2)})}}{\sqrt{2\cdot f_{m,n}}}
>(1+12\,\eta)\sigma \;\Big)}\\
&\le& P\Big(\max_{\stackrel{m_i\le m\le m_{i+1}}{n_j \leq n \leq
n_{j+1}}}(-T_{(m_i,n_j),(m,n)})
>2\eta\,\sigma \sqrt{2\,f_{m_i n_j}}\;\Big)\\
&&\hskip2pc+ P\Big(\max_{\stackrel{m_i\le m\le m_{i+1}}{n_j \leq
n \leq n_{j+1}}}(-T_{(m_i,n_j),(m+a_m^{(1)},n)})
>2\eta\,\sigma \sqrt{2\, f_{m_i n_j} }\;\Big)\\
&&\hskip2pc+ P\Big(\max_{\stackrel{m_i\le m\le m_{i+1}}{n_j \leq
n \leq n_{j+1}}}(-T_{(m_i,n_j),(m,n+a_n^{(2)})})
>2\eta\,\sigma \sqrt{2\,f_{m_i n_j}}\;\Big)\\
&&\hskip2pc+ P\Big(\max_{\stackrel{m_i\le m\le m_{i+1}}{n_j \leq
n \leq n_{j+1}}}T_{(m_i,n_j),(m+a_m^{(1)},n+a_n^{(2)})}
>(1+6\,\eta)\sigma \sqrt{2\,f_{m_i n_j}}\;\Big)\,, \end{eqnarray*} where
$(m_i,n_j)=(e^{\sqrt{ci}},e^{\sqrt{cj}})$ as before. Now, since
$T\stackrel{d}{=}-T$ the second L\'evy inequality for random fields
(\cite{parpar}, Theorem 1) is instrumental and we proceed.
\begin{eqnarray*}
\lefteqn{ P\Big(\max_{\stackrel{m_i\le m\le m_{i+1}}{n_j \leq n
\leq
n_{j+1}}}\frac{T_{(m,n),(m+a_m^{(1)},n+a_n^{(2)})}}{\sqrt{2f_{m,n}}}
>(1+12\,\eta)\sigma \;\Big)}\\
&\le& 2\,P\Big(T_{(m_i,n_j),(m_{i+1},n_{j+1})}
>2\eta \sigma \sqrt{2\,f_{m_i n_j}}\;\Big)\\
&&\hskip2pc+
2\,P\Big(T_{(m_i,n_j),(m_{i+1}+a_{m_{i+1}}^{(1)},n_{j+1})}
>2\eta \sigma \sqrt{2\,f_{m_i n_j}}\;\Big)\\
&&\hskip2pc+
2\,P\Big(T_{(m_i,n_j),(m_{i+1},n_{j+1}+a_{n_{j+1}}^{(2)})}
>2\eta \sigma \sqrt{2\,f_{m_i n_j}} \;\Big)\\
&&\hskip2pc+
2\,P\Big(T_{(m_i,n_j),(m_{i+1}+a_{m_{i+1}}^{(1)},n_{j+1}
+a_{n_{j+1}}^{(2)})}
>(1+6\,\eta)\sigma \sqrt{2\,f_{m_i n_j}}\;\Big)\,, \end{eqnarray*}
By the definition of our subset we have, for all $\eta>0$ and $i,j\ge
i_0(\eta)$ with some integer $i_0\ge 1/c$, that (note $c=\eta^2$),
\bea\label{mnaa}
\begin{cases}
&m_{i+1}-m_i\le  e^{\sqrt{c i}}\,(e^{\sqrt{ci}/(2i)}-1)
\le \eta^2 \, a_{m_i}^{(1)}
\,,\\[4pt]
& n_{j+1}-n_{j}\le e^{\sqrt{c j}}\,(e^{\sqrt{cj}/(2j)}-1)
\le\eta^2\,a_{n_j}^{(2)}\,,\\[4pt]
&a_{m_{i+1}}^{(1)}/a_{m_i}^{(1)}\le \exp{\sqrt{c/i}}\le 1+2\eta^2\,,\\[4pt]
&a_{n_{j+1}}^{(2)}/a_{n_j}^{(2)} \le \exp{\sqrt{c/j}}\le
1+2\eta^2\,,\end{cases}
\eea
from which it follows that the variances satisfy
\bea\label{varianslog}\begin{cases}
&\var(T_{(m_i,n_j),(m_{i+1},n_{j+1} )})\le\sigma^2\,\eta^4
a_{m_i}^{(1)}\,a_{n_j}^{(2)}=\sigma^2\,\eta^4 a_{m_i,n_j} \,,\\[3pt]
&\var(T_{(m_i,n_j),(m_{i+1}+a_{m_{i+1}}^{(1)},n_{j+1})})
\le\sigma^2(1+3\eta^2)\,\eta^2
a_{m_i}^{(1)}\,a_{n_j}^{(2)}=\sigma^2(1+3\eta^2)\,\eta^2 a_{m_i,n_j} \,,\\[3pt]
&\var(T_{(m_i,n_j),(m_{i+1},n_{j+1}+a_{n_{j+1}}^{(2)})})
\le\sigma^2(1+3\eta^2)\,\eta^2
a_{m_i}^{(1)}\,a_{n_j}^{(2)}=\sigma^2(1+3\eta^2)\,\eta^2 a_{m_i,n_j} \,,\\[3pt]
&\var(T_{(m_i,n_j),(m_{i+1}+a_{m_{i+1}}^{(1)},n_{j+1}
+a_{n_{j+1}}^{(2)})})\le\sigma^2(1+3\eta^2)^2
a_{m_i}^{(1)}\,a_{n_j}^{(2)}= \sigma^2(1+3\eta^2)^2 a_{m_i,n_j}\,, \end{cases}
\eea
for  $i,\,j\ge i_0$.

From here on the procedure from the previous sections
applies for these terms. For example, the second last probability can be
bounded with respect to $T'$ by
\beaa&&\hskip-2pc 2\,P\Big(T_{(m_i,n_j),(m_i + \eta^4
a^{(1)}_{m_i},n_j + (1+3\eta^2)a_{n_{j}}^{(2)})}'
>2\eta \sigma {\sqrt{2f_{m_i,n_j}}}\Big) \\[1mm]
&&\hskip2pc\leq
\exp\Big\{-\frac{(1-\delta)^3 4\eta^2\sigma^2 2\log(\log m_i\log n_j)}
{\sigma^2 \eta^2(1+3\eta^2)}\Big\}
\le C (ij)^{-4(1-\delta)^3}\,,
\eeaa
provided $\eta^2<1/3$, and  the last term by
\beaa&&\hskip-2pc 2\,P\Big(T_{(m_i,n_j),(m_i +
(1+2\eta^2) a^{(1)}_{m_i},n_j + (1+2\eta^2)a_{n_{j}}^{(2)})}'
>2\eta \sigma {\sqrt{2f_{m_i,n_j}}}\Big) \\[1mm]
&&\hskip2pc\leq \exp\Big\{-\frac{(1-\delta)^3 (1+6\eta)^2 \sigma^2 2\log(\log
m_i\log n_j)}{\sigma^2 (1+3\eta^2)^2}\Big\}
\le C(i j)^{-((1+6\eta)/(1+3\eta^2))^2(1-\delta)^3}\,.
\eeaa
Both cases yield summable double sequences if $\delta$ and $\eta$ are
small enough.

Since constants are not relevant for $T'',T'''$ we finally obtain
\[
\sum_{k,j\ge i_0}P\Big(\max_{\stackrel{m_i\le m\le
m_{i+1}}{n_j \leq n \leq n_{j+1}}}
\frac{T_{(m,n),(m+a_m^{(1)},n+a_n^{(2)})}}{\sqrt{2f_{m,n}}}
>(1+12\,\eta)\sigma \Big)
<\infty\,,
\]
which implies that the $\limsup \dots \le (1+12\eta)\sigma$, and
since $\eta$ can be chosen arbitrarily small the upper bound is
proved for the entire field.

Desymmetrization follows along the
usual arguments; note that
\[ \frac{E\left|T_{(m,n),(m+a_m^{(1)},n+a_n^{(2)})}\right|}{\sqrt{f_{m,n}}}
\leq \sqrt{\frac{E (T_{(m,n),(m+a_m^{(1)},n+a_n^{(2)})})^2}{f_{m,n}}}=o(1)
\ttt{as}m,n \to \infty\,.\]

\subsubsection*{Lower bound}

Let $c>2$ and define $m_i=n_i=\exp{\{\sqrt{ci}\}}$. We first note that
\[\frac{e^{\sqrt{ci}}+e^{\sqrt{ci}}/\sqrt{ci}}{e^{\sqrt{c(i+1)}}}
= \frac{1+1/\sqrt{ci}}{\exp{\big(\sqrt{ci}(\sqrt{1+1/i}-1)\big)}}\sim
\frac{1+1/\sqrt{ci}}{1+\sqrt{c/4i}} <1\,,\]
eventually and, hence, that the windows of the subset
$(m_i,\,n_j)$, for $i,j$ large, are disjoint, which
means that different blocks (eventually)
consist of independent random variables. Using
the lower exponential bound in (\ref{expbdslog}) we then obtain
\beaa
&&\hskip-3pc P\Big(T_{(m_i,n_j)\,,\,(m_i+m_i/\log m_i,n_j+n_j/\log n_j)}'
>(1-\eta)\,\sigma\sqrt{2 f_{m_i,n_j}}\;\Big)\\
&&\hskip2pc \geq \exp \Big\{\frac{(1-\eta)^2
(1+\delta)^2(1+\gamma)}{(1-\delta)}\cdot(2 \log\log m_i +2 \log\log n_j)
\Big\}\\
&&\hskip2pc \geq
C\,(ij)^{-\frac{(1-\eta)^2\,(1+\delta)^2\,(1+\gamma)}{(1-\delta)}}\,,
\eeaa
which is a divergent minorant for a choice of $\eta>0 $ such that $1-\eta<
\sqrt{\frac{1-\delta}{(1+\delta)^2(1+\gamma)}}$. Now, since
$\delta$ and $\gamma$ can be chosen arbitrarily small the same is
true for $\eta$. The desired result finally follows via the Borel--Cantelli
Lemma and the fact that $T''$ and $T'''$ are small.

Strictly speaking, this provides the desired lower bound for our
\emph{subset}, after which the overall lower bound follows from
that fact that the limsup over \emph{all\/} windows is at least as
large as the limsup over a subset.

\subsubsection*{Necessity}

It follows from the assumption  that
\[\limsup_{m,n\to\infty}\frac{|X_{m,n}|}{\sqrt{f_{m,n}}}=
\limsup_{m,n\to\infty}\frac{|X_{m,n}|}{\sqrt{\frac{2mn(\log\log m+\log\log n)}
{\log m\log n}}}<\infty\ttt{a.s.,}\] from which the necessity of the
moment assumption is immediate in view of (\ref{corrlog}). An
application of the sufficiency and the strong law of large numbers
then tells us that the mean must be equal to
zero.

\subsection{Proof of Theorem \ref{thmmain}(ii)}
\subsubsection*{Sufficiency}
The proof follows the same procedure with obvious modifications. We have
\[a_{m,n}=\frac{mn}{\log\log m\log\log n},\qquad
d_{m,n}= \log\log m + \log\log n,\qquad f_{m,n}= \frac{mn(\log\log m+
\log\log n)}{\log\log m \log\log n}\,,\] and the Kolmogorov exponential bounds
(\ref{expbds}) yield
\beaa&&\hskip-4pc
P(T_{(m,n)\,,\,(m+m/\log\log m,n+n/\log\log n}'
>\varepsilon\sqrt{2mn\frac{\log\log m +\log\log n}{\log\log m \log\log n}}\\
&&\hskip2pc\begin{cases}&\!\!
\leq \exp\big\{-\frac{\varepsilon^2(1-\delta)^3}{\sigma^2}
\cdot (\log\log m+\log\log n)\;\big\}\,,\\[6pt]
&\!\!\geq \exp\big\{-\frac{\varepsilon^2(1+\delta)^2(1+\gamma)}
{\sigma^2(1-\delta)}\cdot (\log\log m +\log\log n)\big\}\,,
\end{cases}
\eeaa
for $n$ large and any $\gamma>0$.

For the upper bound we define, with some constant $c>0$ to be defined
later,
\[(m_i,n_j)=\big(e^{c\,i/\log(i+1)},e^{c\,j/\log (j+1)}\big)
\ttt{for}i,j \ge \max\{\log (1/c)/c,1\}\,,\]
to obtain
\beaa &&\hskip-3pc P\Big(T_{(m_i,n_j)\,,\,(m_i+m_i/\log\log
m_i,n_j+n_j/\log\log n_j)}'
>\varepsilon\sqrt{2m_i n_j\frac{\log\log m_i +\log\log n_j}
{\log\log m_i \log\log n_j}}\;\Big)\\
&&\hskip3pc \leq C\,\left(\frac{\log i\log
j}{ij}\right)^{\,-\,\frac{\varepsilon^2(1-\delta)^3}{\sigma^2}}\,.
\eeaa
The analog of (\ref{pzweiprel}) turns into
\[
P''(m,n)\leq C\bigg(\frac{\log\log m+\log\log n}{\log (mn)}\bigg)^N,
\]
so that, along the current subset, we obtain
\[
P''(m_i,n_j)\leq C\bigg(\frac{(\log i+\log j)\log i\,\log
j}{i+j}\bigg)^N\,,
\]
and therefore, in complete analogy with the previous case (use $N\ge 3$ here),
\[
\limsup_{i,j\to\infty}
\frac{|T_{(m_i,n_j)\,,\,(m_i+m_i/\log\log m_i,n_j+n_j/\log\log n_j)}''|}
{\sqrt{f_{m_i,n_j}}}\leq \delta \ttt{a.s.}
\]
The argument for $T'''$ is the same as before, so that, recall
(\ref{corrlog}),
\[
\sum_{m,n}P(|X_{m,n}|>\eta\sqrt{f_{m,n}}) =\sum_{m,n}
P\Big(|X|>\eta\sqrt{2mn\frac{\log\log m +\log\log n}{\log\log m
\log\log n}}\;\Big)<\infty\,,
\]
since, by assumption, $E\,X^2\log^+|X|\log^+\log^+ |X|<\infty$.

Finally, by combining the pieces, it follows that
\bea\label{zloglog} \limsup_{i,j\to\infty}
\frac{T_{(m_i,n_j)\,,\,(m_i+m_i/\log\log m_i,n_j+n_j/\log\log n_j)}}
{\sqrt{2f_{m_i,n_j}}} \leq  \sigma(1-\delta)^{-3/2}+\delta
\ttt{a.s.,} \eea for
$(m_i,n_j)=\big(e^{c\,i/\log(i)},\,e^{c\,j/\log(j)}\big)$ with any
$c>0$.

\subsubsection*{Sufficiency for the entire field}
Once again we have proved the theorem for a suitable subset, and
it remains to show that our process behaves accordingly for the entire
field. In order to achieve this we use the same procedure as
for the corresponding part of the proof of part (i) with
logarithms replaced by iterated logarithms.

We first note that here,
\begin{eqnarray*}
\frac{m_{i+1}-m_i}{a_{m_i}}&\le&
\frac{\exp{\left(\frac{c(i+1)}{\log (i+1)}\right)}-\exp{\left(\frac{ci}
{\log (i+1)}\right)}}
{\exp{\left(\frac{ci}{\log (i+1)}\right)}}\log(ci/\log(i+1))\le 2c\,,\\
 \frac{a_{m_{i+1}}}{a_{m_i}}
&\le& \left(\frac{\exp{\left(\frac{c(i+1)}{\log (i+1)}\right)}}
{\exp{\left(\frac{ci}{\log (i+1)}\right)}}\right)\,
\frac{\log(ci/\log(i+1))}{\log(c(i+1)/\log(i+2))}\le 1+2c\,,
\end{eqnarray*}
for all integers $i$ considered, so that, by choosing $c=\eta^2/2$,
the analog of (\ref{varianslog}) becomes
\bea\label{varianloglog}\begin{cases}
&\var(T_{(m_i,n_j),(m_{i+1},n_{j+1})})\le\sigma^2\,\eta^4 a_{m_i,n_j}\,,\\[3pt]
&\var(T_{(m_i,n_j),(m_{i+1}+a_{m_{i+1}}^{(1)},n_{j+1})})
\le\sigma^2(1+2\eta^2)\,\eta^2 a_{m_i,n_j} \,,\\[3pt]
&\var(T_{(m_i,n_j),(m_{i+1},n_{j+1}+a_{n_{j+1}}^{(2)})})
\le \sigma^2(1+2\eta^2)\,\eta^2 a_{m_i,n_j} \,,\\[3pt]
&\var(T_{(m_i,n_j),(m_{i+1}+a_{m_{i+1}}^{(1)},n_{j+1}
+a_{n_{j+1}}^{(2)})})\le\sigma^2(1+2\eta^2)^2 a_{m_i,n_j}\,, \end{cases}
\eea
for  $i,\,j\ge i_0$.

We omit the remaining details.

\subsubsection*{Lower bound}

Let $c>1$ and define $m_i=n_i=\exp{\{ci/\log (i+1)\}}$. Now,
\beaa
&&\hskip-2pc\frac{e^{ci/\log (i+1)}+e^{ci/\log(i+1)}/\log
(ci/\log(i+1))} {e^{c(i+1)/\log((i+2)}}=
\frac{1+1/\log(ci/\log(i+1))}
{\exp{\big(c(i+1)/\log\{(i+1)(1+1/(i+1))\}\big)}}\\
&&\hskip2pc\sim
\frac{1+1/\log(ci)}{\exp{\big(c/\log(i+1)\big)}}<1\,, \eeaa
eventually, which, in analogy to above tells us that the windows
along the subset $(m_i,\,n_j)$, for $i,j$ large, are disjoint, that is, that
different blocks consist of independent random variables. The lower
exponential bound in (\ref{expbdslog}) then takes care of $T'$,
after which the conclusion of this part of the proof is analogous to
the conclusion of the lower bound for part (i) of the theorem.

\subsubsection*{Necessity}
In this case
\[\limsup_{m,n\to\infty}\frac{|X_{m,n}|}{\sqrt{f_{m,n}}}=
\limsup_{m,n\to\infty}\frac{|X_{m,n}|}{\sqrt{mn\frac{\log\log m +\log\log n}
{\log\log m \log\log n}}}<\infty\ttt{a.s.,}\]
which, in view of (\ref{corrloglog}), is equivalent $E\,X^2\log^+|X|
\log^+\log^+|X| <\infty$. The remaining part follows as before.

\subsection{Proof of Theorem \ref{thmmain}(iii)}
\subsubsection*{Sufficiency}

The proof follows the same pattern as that of the previous ones,
although the rates of the stretches of the windows now differ in the
two directions.  Accordingly, the components of the subset will grow
at different rates. We have
\[a_{m,n}=\frac{mn}{\log m\log\log n},\qquad
d_{m,n}= 2( \log\log m+\log\log n),\qquad f_{m,n}=
\frac{2mn\,(\log\log m +\log\log n)}{\log m\,\log\log n}\,.\]
The Kolmogorov exponential upper bound along the subset
$(m_i,n_j)=\big(e^{\sqrt{ci}},e^{\sqrt{cj}}\big)$, $i,j\geq\max\{1/c\,,1\}$
with an arbitrarily small $c>0$, then tells us that
\[
P\Big(T_{(m_i,n_j)\,,\,(m_i+m_i/\log m_i,n_j+n_j/\log\log n_j)}'
>\varepsilon\sqrt{2 f_{m_i,n_j}}\Big)
\leq C (ij)^{\,-\,\frac{\varepsilon^2(1-\delta)^3}{\sigma^2}}\,.
\]
As for $P''$ we obtain
\[
P''(m_i,n_j)\leq C\bigg(\frac{\log\sqrt{i}+ \log
\sqrt{j}}{(\sqrt{i}+\sqrt{j})^2} \bigg)^N\,,
\]
and, hence, that for $N\ge 3$,
\[
\limsup_{i,j\to\infty}
\frac{|T_{(m_i,n_j)\,,\,(m_i+m_i/\log m_i,n_j+n_j/\log\log n_j)}''|}
{\sqrt{f_{m_i,n_j}}}\leq \delta \ttt{a.s.,}
\]
and for $T'''$, via (\ref{corrlogloglog}) and the fact that
$E\,X^2\,(\log^+|X|)^2<\infty$,
\[
\sum_{m,n}P(|X_{m,n}|>\eta\sqrt{f_{m,n}}) =\sum_{m,n}
P\Big(|X|>\eta\sqrt{\frac{mn(\log\log m+\log\log n)}{\log m
\log\log n}}\;\Big)<\infty\,.
\]
Combining everything finally yields
\bea\label{zlogloglog} \limsup_{i,j\to\infty}
\frac{T_{(m_i,n_j)\,,\,(m_i+m_i/\log m_i,n_j+n_j/\log\log n_j)}}
{\sqrt{2f_{m_i,n_j}}}
\leq  \sigma(1-\delta)^{-3/2}+\delta \ttt{a.s.}
\eea

\subsubsection*{Sufficiency for the entire field}
This case is similar to the first one. Since the subsequence for the
$n$-coordinate is denser than in the $\log\log$-case in part (ii),
the task to fill the gaps is even simpler than before.

\subsubsection*{Lower bound}

Once again we let $c>2$ and set $m_i=n_i=e^{\sqrt{ci}}$. For the first
coordinate we have $m_i+m_i/\log m_i< m_{i+1}$ eventually, i.e.,
the blocks consist of independent random variables. The lower
exponential bound now yields
\[\sum_{i,j} P\Big(T_{(m_i,n_j),(m_i+m_i/\log m_i,n_j+n_j/\log\log
n_j)}'> (1-\eta)\sigma \sqrt{4\frac{mn(\log\log m_i +\log\log
n_j)}{\log m_i\log\log n_j}}\;\Big)=\infty\,,\] and the lower
bound follows as before.

\subsubsection*{Necessity}
This time
\[\limsup_{m,n\to\infty}\frac{|X_{m,n}|}{\sqrt{f_{m,n}}}
=\limsup_{m,n\to\infty}\frac{|X_{m,n}|}
{\sqrt{\frac{2mn(\log\log m +\log\log n)}
{\log m \log\log n}}}<\infty\ttt{a.s.,}
\]
and the necessity follows along the usual lines.

\section{Additional results}
\setcounter{equation}{0}\label{7075}

In this section we consider windows of the form
\[T_{(m,n)\,,\,(m+m^{\alpha},\,n+n/\log n)},\]
that is, the windows are rectangles where one side has a length
$a_n^{(2)}=n/\log n$ as before, whereas the other one has length
$a_n^{(1)}=m^\alpha$ for some $\alpha\in(0,1)$, and, hence, is
\emph{much\/} shorter, so that, noticing that $r(n)=(1-\alpha)\log n$,
our usual quantities become
\[a_{m,n}=\frac{m^{\alpha}n}{\log n},\qquad
d_{m,n}=(1-\alpha)(\log m + \log n),\qquad
f_{m,n}=(1-\alpha)\frac{m^{\alpha}n(\log m+\log n)}{\log n}.\]
Here is now our result for this setting.
\begin{theorem} Let $0<\alpha< 1$. If
\[E\,X^{2/\alpha}(\log^+|X|)^{-1/\alpha}<\infty,\ttt{and} E(X)=0,\;
E(X^2)=\sigma^2\,,
\]
then
\[\limsup_{m,n\to\infty}
\frac{T_{(m,n)\,,\,(m+m^{\alpha},\,n+n/\log n)}}
{\sqrt{2 m^\alpha n\frac{(1-\alpha)\log(mn)}{\log n}}}=\sigma\quad\mbox{a.s.}
\]
Conversely, if
\[P\Big(\limsup_{m,n\to\infty}
\frac{|T_{(m,n)\,,\,(m+m^{\alpha},\,n+n/\log n)}|}
{\sqrt{ m^\alpha n\frac{\log(mn)}{\log n}}}<\infty\Big)>0\,,\]
then $E\,X^{2/\alpha}(\log^+|X|)^{-1/\alpha}<\infty$ and $E(X)=0$.
\end{theorem}
\begin{proof} As for the moment condition,
\bea\label{nalpha}
\limsup_{m,n\to\infty}\frac{|X_{m,n}|}{\sqrt{m^\alpha n\frac{\log(mn)}
{\log n}}}<\infty\quad\mbox{a.s.}
&\iff&\sum_{m,n}P\Big(|X|>\sqrt{m^\alpha n\frac{\log(mn)}
{\log n}}\Big)<\infty\nonumber\\
&\iff& E\,(X^2/\log^+ |X|)^{1/\alpha}<\infty\,,
\eea
which will be verified in the Appendix.

The next step is to truncate and split the window into the usual three parts.
Toward that end we note that (\ref{bn}) becomes
\bea\label{bn7075} b_{m,n}
  =\frac{\sigma\delta}{\varepsilon} \sqrt{\frac{m^\alpha n}
{(1-\alpha)\big(\log m\,\log n+(\log n)^2\big)}}\,,\eea
after which we introduce the standard three components
$X_{m,n}'$, $X_{m,n}''$, and $X_{m,n}'''$ of the summand $X_{m,n}$
and observe that
\[\frac{a_{m,n}}{b_{m,n}}=\frac{\sigma\delta}{\varepsilon}\sqrt{f_{m,n}}\,,\]
so that, via the usual estimates,
\[|E\,T_{(m,n)\,,\,(m+m^{\alpha},\,n+n/\log n)}'|=o(\sqrt{f_{m,n}})
\ttt{as}m,n\to\infty.\]
With this in mind the exponential bounds yield
\bea
P\Big(T_{(m,n)\,,\,(m+m^{\alpha},\,n+n/\log n)}'
>\varepsilon\sqrt{2f_{m,n}}\Big)
\;\begin{cases}&
\leq \exp\big\{-\frac{\varepsilon^2(1-\delta)^3}{\sigma^2}
\cdot (1-\alpha)\log(mn)\big\}
\,,\\[6pt]
&\geq \exp\big\{-\frac{\varepsilon^2(1+\delta)^2(1+\gamma)}{\sigma^2(1-\delta)}
\cdot (1-\alpha)\log(mn)\big\}\,,
\end{cases}\label{expbdsalpha}
\eea
for any $\gamma>0$ and $\max\{m,n\}$ large.

We note, in passing, that the sum of the probabilities converges for
suitably chosen geometrically increasing sequences.

As for $T''$, (\ref{ttzwei}) turns into
\beaa
P(|T_{(m,n)\,,\,(m+m^{\alpha},\,n+n/\log n)}''|
>\eta\sqrt{f_{m,n}})&\leq&C\Big(\frac{m^{\alpha} n}{\log n}\Big)^N
\bigg(\frac{E|X|^{2/\alpha}}{\big(\frac{m^{\alpha}n}
{(\log n)^2+\log m\,\log n}\big)^{1/\alpha}}\bigg)^N\\[4pt]
&=&C\frac{\big((\log n)^{1-\alpha}(\log n+\log m)\big)^{N/\alpha}}
{m^{(1-\alpha)N} n^{((1/\alpha)-1)N}}\,,
\eeaa
so that
\bea\label{pzweialpha}
\sum_{m,n}P(|T_{(m,n)\,,\,(m+m^{\alpha},\,n+n/\log n)}''|
>\eta\sqrt{f_{m,n}})<\infty
\eea
whenever $N>\max\{1/(1-\alpha),\alpha/(1-\alpha)\}=1/(1-\alpha)$, that
is (cf.\ \cite{gs1}), recalling that $N\delta>\eta$ as before,
\[\limsup_{m,n\to\infty}\frac{|T_{(m,n)\,,\,(m+m^{\alpha},\,n+n/\log n)}''|}
{\sqrt{f_{m,n}}}\leq\frac{\delta}{1-\alpha}.\]
By copying the arguments for $T'''$ from above (as well as from our
predecessors) it finally follows that
\[\frac{T_{(m,n)\,,\,(m+m^{\alpha},\,n+n/\log n)}'''}
{\sqrt{f_{m,n}}}\asto0\ttt{as}m,n\to\infty.
\]
\subsubsection*{Upper bound}
For the upper bound, let $c>0$ and consider the subset
\bea\label{subupp} (m_i, n_j)=\Big(c\,i^{1/(1-\alpha)},\, c\,
j^{1/(1-\alpha)}\Big) \ttt{for}i,j\geq \max\{c^{\alpha-1}\,,1\}. \eea
Summarizing the three contributions, and taking the arbitrariness
of $\delta$ and $\eta$ into account, we have shown that for this
subset we have \bea\label{limsupeins}\limsup_{i,j\to\infty}
\frac{T_{(m_i,n_j)\,,\,(m_i+m_i^{\alpha},\,n_j+n_j/\log n_j)}}
{\sqrt{2f_{m_i,n_j}}}\leq\sigma\ttt{a.s.} \eea
for any $c>0$.

\subsubsection*{Sufficiency for the entire field}
Here we have that for $i,j\ge c^{2(\alpha-1)}$ and
$\log j/j\leq 1/\sqrt{j}$ with $c>0$, sufficiently small,

\bea\label{mnmix}
 \begin{cases}
&m_{i+1}-m_i\sim\dfrac{c}{1-\alpha}i^{\alpha/(1-\alpha)}
=\dfrac{c^{1-\alpha}}{1-\alpha}\,a_{m_i}^{(1)}\,,\\[8pt]
& n_{j+1}-n_{j}\sim\dfrac{c}{1-\alpha}j^{\alpha/(1-\alpha)}
=\dfrac{\log j}{(1-\alpha)^2j}\,a_{n_j}^{(2)}
\leq \dfrac{c^{1-\alpha}}{(1-\alpha)^2} \,a_{n_j}^{(2)} \,,\\[8pt]
&a_{m_{i+1}}^{(1)}/a_{m_i}^{(1)}
\sim 1+\dfrac{\alpha}{1-\alpha}\cdot\dfrac1{i}
\leq 1+\dfrac{\alpha c^{2(1-\alpha)}}{1-\alpha}\,,\\[8pt]
&a_{n_{j+1}}^{(2)}/a_{n_j}^{(2)}=\Big(1+\dfrac1{j}\Big)^{1/(1-\alpha)}
\cdot\dfrac{\log j}{\log(j+1)} \leq1+\dfrac1{1-\alpha}\cdot\dfrac1{j}
\leq 1+\dfrac{c^{2(1-\alpha)}}{1-\alpha}\,.
\end{cases}
\eea
By choosing $\eta^2=c^{1-\alpha}/(1-\alpha)^2$ (which implies that
$c=((1-\alpha)^2\eta^2)^{1/(1-\alpha)}\leq \eta^2$), the estimates for the
variances turn out as follows:
\bea\label{variansmix}\begin{cases}
&\var(T_{(m_i,n_j),(m_{i+1},n_{j+1} )})\le\sigma^2\,
\dfrac{c^{2(1-\alpha)}}{(1-\alpha)^3}\,
a_{m_i,n_j}\le\sigma^2\,
\eta^4\,a_{m_i,n_j}  \,,\\[8pt]
&\var(T_{(m_i,n_j),(m_{i+1}+a_{m_{i+1}}^{(1)},n_{j+1})})
\le\sigma^2\Big(1+\dfrac{c^{1-\alpha}+\alpha c^{2(1-\alpha)}}{1-\alpha}\Big)
\dfrac{c^{1-\alpha}}{(1-\alpha)^2}\, a_{m_i,n_j}
\\[6pt]
&\hskip4.7cm\leq\sigma^2\Big(1+\dfrac{2\,c^{1-\alpha}}{(1-\alpha)^2}\Big)
\dfrac{c^{1-\alpha}}{(1-\alpha)^2}\, a_{m_i,n_j}\\[6pt]
&\hskip4.7cm\leq\sigma^2 (1+\eta^2)\,\eta^2\, a_{m_i,n_j} \,,\\[8pt]
&\var(T_{(m_i,n_j),(m_{i+1},n_{j+1}+a_{n_{j+1}}^{(2)})})
\le\sigma^2\dfrac{c^{1-\alpha}}{1-\alpha}\Big(
\dfrac{c^{1-\alpha}}{(1-\alpha)^2}+1+\dfrac{c^{2(1-\alpha)}}{1-\alpha}\Big)\,
a_{m_i,n_j} \\[6pt]
&\hskip4.7cm\leq  \sigma^2\dfrac{c^{1-\alpha}}{1-\alpha}
\Big(1+\dfrac{2\,c^{1-\alpha}}{(1-\alpha)^2}\Big)a_{m_i,n_j} \\[6pt]
&\hskip4.7cm\leq  \sigma^2 \eta^2\big(1+2\eta^2\big)\,a_{m_i,n_j}\,,\\[8pt]
&\var(T_{(m_i,n_j),(m_{i+1}+a_{m_{i+1}}^{(1)},n_{j+1}
+a_{n_{j+1}}^{(2)})})
\le\sigma^2\Big(1+\dfrac{c^{1-\alpha}+\alpha c^{2(1-\alpha)}}{1-\alpha}\Big)
\\[6pt]
&\hskip7cm \times\Big(
\dfrac{c^{1-\alpha}}{(1-\alpha)^2}+1+\dfrac{c^{2(1-\alpha)}}{1-\alpha}\Big)
 a_{m_i,n_j}\\[6pt]
&\hskip5.8cm
\leq\sigma^2\Big(1+\dfrac{2\,c^{1-\alpha}}{(1-\alpha)^2}\Big)^2
a_{m_i,n_j}\\[6pt]
&\hskip5.8cm
\leq\sigma^2 \big(1+2\eta^2\big)^2
a_{m_i,n_j}\,, \end{cases}
\eea
for  $i,\,j\ge i_0$.

The remaining part of the proof follows the usual procedure with obvious
changes and is omitted.

\subsubsection*{Lower bound}
For the lower bound we set, with some $c>1$,
\[m_i=n_i=c\,i^{1/(1-\alpha)}\,,\]
so that combining the contributions above, taking the arbitrariness of
$\delta$, $\gamma$, and $\eta$ into account, we obtain
\bea\label{limsupzwei}\limsup_{i,j\to\infty}
\frac{T_{(m_i,n_j)\,,\,(m_i+m_i^{\alpha},\,n_j+n_j/\log n_j)}}
{\sqrt{2f_{m_i,n_j}}}\geq\sigma\ttt{a.s.}
\eea
Namely, since, for the current subset we have
$m_i+m_i^{\alpha}<m_{i+1}$ eventually, the windows of the subset are
disjoint, so that (\ref{limsupzwei}) follows in view of the second
Borel--Cantelli lemma, after which, once again, the fact that
the limsup over \emph{all\/}
windows is at least as large as the limsup over a subset establishes
(\ref{limsupzwei}).

\subsubsection*{Necessity}
This part follows the standard pattern.
\end{proof}

\section{Some final remarks}
\setcounter{equation}{0}\label{anm}

\begin{enumerate} \item The subsequences we have used in the proof of the
upper bound have, throughout been the same for both coordinates, even
in those cases when the edges of the windows grow at different rates,
viz., in the proofs of Theorem \ref{thmmain}(iii) and
(\ref{subupp}). Intuitively one might imagine subsequences that
grow at different rates to be more natural. However, since there is less
cancellation going on along the direction of the shorter edge, that
is, since the fluctuations along the direction of the shorter edge are
stronger, it turns out that the shortest edge
determines the moment condition (as in \cite{gs2}). Our choice of
subset implies that the windows of the subset are ``covered'' many
times along the longer direction, which, on the one hand is not
necessary, but on the other does not exhibit any ``harmful waste'', in
that our estimates produce the best result.

\item Analogous to our previous results \cite{gs1,gs2} we can obtain limit
theorems over subsets, which, in particular, show that all reals
between the limit inferior and the limit superior are limit points.

\item An analysis of the proofs shows that it is also possible to
formulate a more general
result. However, the proof is more tedious and the moment condition is not so
explicit.\\\phantom{osos}
Assume that
$L_i(x)\nearrow \infty$ for $ i=1,2$ are slowly varying
functions satisfying (\ref{an1}), and suppose that $L_1(x)\ge L_2(x)$
for $x> \mbox{some } x_0$. Further, define $r(n)=\log(L_1(n)\,\log n)$ and
\[G(u,v)=\frac{u\,v}{L_1(u)\,L_2(v)}(r(u)+r(v))\]
with the associated $M$-function
\[
M(x)\,=\,\big|\{(u,v)\in [1,\,\infty)
\times [1,\,\infty):G(u,v)\le x\}\big|.
\]
Then the following result holds:
\begin{theorem}  If $E\,M(X^2)<\infty$, $E(X)=0$, and $E(X^2)=\sigma^2$, then
\bea\label{allmsats}
\limsup_{m,n\to\infty}\, \frac{T_{(m,n)\,,\,(m+m/L_1(m),\, n+n/L_2(n))}}
{\sqrt{2 G(m,n)}}=\sigma\ttt{a.s.}\eea
Conversely, if the $\limsup$ is finite with positive probability,
then $E\,M(X^2)<\infty$, $E(X)=0$, and (\ref{allmsats}) holds with
$\var X=\sigma^2$.
\end{theorem}
\begin{remark} \emph{Note that $M(x)/x \to \infty$, that is, the second moment
exists if $E\,M(X^2)<\infty$.}\end{remark}
\begin{remark} \emph{
Due to the generality it is not possible to
describe  the function $M(x)$ more
precisely, but one can prove (see the Appendix) that
\begin{equation}\label{UnglMG}
M(x)\le C\,\frac{x\,L_2(x)}{r(x)}\int^{x\,L_1(x)/r(x)}\frac{L_1(u)}{u}\,du\,.
\end{equation}
The proof follows the lines of the proofs above using the techniques
from \cite{gjs} with subsequences $m_k=n_k=\psi(ck)$ where
$\psi(x)=\varphi^{-1}(x)$ with $\varphi(x)=\int^x L_1(u)/u\, du
\nearrow \infty$ and being in $\cal{SV}$.}\end{remark}
\begin{remark}\vspace{-4mm} \emph{For Theorem \ref{thmmain} the estimate
in (\ref{UnglMG}) is precise enough to give the correct moment
assumptions.}\vsb
\end{remark}
\end{enumerate}

\renewcommand{\theequation}{\Alph{section}.\arabic{equation}}
\appendix
\section*{Appendix}
\setcounter{equation}{0}
\section*{Proofs of (\ref{corrlog}) -- (\ref{corrlogloglog}), (\ref{nalpha})
and (\ref{UnglMG})}
In this appendix we present the ``elementary but tedious
calculations'' needed for the verification of the second half of
Corollary \ref{corrmom} and the further analogs.
Throughout it is tacitly assumed that all logarithms and
iterated logarithms are to be interpreted as $\max\{\log x,1\}$ and
$\max\{\log\log x,1\}$, respectively. Moreover, we shall use
$a(x)\approx b(x)$ to denote that $0< c_1\le a(x)/b(x)\le c_2$.
Finally, the constant $C$ below may be different at each occurrence.

\subsubsection*{Proof of  (\ref{corrlog})}

Since $\log u/u\searrow$ we obtain \beaa I_1(u)&=& \iint_{xy(\log\log
x+\log\log y)/\log x\log y\leq u}1\,dxdy
=\iint_{xy(\log(\log x\log y))/\log x\log y\leq u}1\,dxdy\\[4pt]
&\leq&\iint_{xy\leq u(\log u)^2/\log((\log u)^2)}1\,dxdy
\le C \iint_{xy\leq u(\log u)^2/\log\log u}1\,dxdy
\le C u(\log u)^3/\log\log u, \eeaa
and
\beaa
I_1(u)&\geq&\int_{\sqrt{u}}^u \int_{x \geq y^{-1}
u(\log\sqrt{u})^2/\log((\log \sqrt{u})^2)}1\,dxdy
\geq C\int_{\sqrt{u}}^u\int_{x\leq (4y)^{-1}u(\log u)^2/\log\log u}1\,dxdy\\
&\geq& C u(\log u)^3/\log\log u.
\eeaa
This shows that $I_1(u)\approx\, u(\log u)^3/\log\log u$ as
$u\to\infty$, i.e., we may choose $M(x)=x\,(\log x)^3/\log\log
x$, and (\ref{corrlog}) follows.

\subsubsection*{Proof of  (\ref{corrloglog})}
We have \beaa I_2(u)&=&
\iint_{xy\log(\log x\cdot\log y)/\log\log x\log\log y)\leq u}1\,dxdy
=\iint_{xy(\frac1{\log\log x}+\frac1{\log\log y})\leq u}1\,dxdy\\[4pt]
&\leq&\iint_{xy\leq 2u\log\log u}1\,dxdy \leq C u\log u\log\log u, \eeaa
and
\[
I_2(u)\geq\int_{\sqrt{u}}^u \int_{x \leq 2\,y^{-1}
u\log\log\sqrt{u}}1\,dxdy \geq C u\log u\log\log u\, .
\]
This shows that $I_2(u)\approx \log u\log\log u$ as $u\to\infty$,
from which we conclude that (\ref{corrloglog}) holds.

\subsubsection*{Proof of  (\ref{corrlogloglog})}

We first note that, $\frac{\log\log x}{\log\log y}+1\geq1$, and that
for $\sqrt{u}\leq x,y\leq u$,
\[\frac{\log\log x}{\log\log y}+1\leq\frac{\log\log u}{\log\log \sqrt{u}}+1
=\frac{\log\log u}{- \log 2+ \log\log u }+1 =2+\frac{\log 2}{-\log
2+\log\log u}\leq C,\] so that, on the one hand, \beaa
I_3(u)&=&\iint_{xy(\log\log x+\log\log y)/\log x\log\log y)\leq
u}1\,dxdy = \iint_{\frac{xy}{\log x}(\frac{\log\log x}{\log\log
y}+1)\leq u}1\,dxdy
\\[4pt]
&\leq&\iint_{\frac{xy}{\log x}\leq u}1\,dxdy
\leq \int_1^{u \log u} u\, \frac{\log x}{x}\,dx
= \frac{u(\log u +\log\log u)^2}{2}\,,
\eeaa
and on the other,
\[
I_3(u)\geq\iint_{\frac{Cxy}{\log x}\leq u}1\,dxdy\sim C u(\log u)^2\,.
\]
This shows that $I_3(u)\approx u(\log u)^2$ as $u\to\infty$, which
proves (\ref{corrlogloglog}).

\subsubsection*{Proof of (\ref{nalpha})}
We first note that, since $1<x,y<u$, we have
\[x^{\alpha}y\frac{\log x}{\log y}\leq
x^{\alpha}y\Big(\frac{\log x}{\log y}+1\Big)
=x^{\alpha}y\frac{\log x+\log y}{\log y}
\leq 2x^{\alpha}y\frac{\log u}{\log y}\,.\]
Thus,
\beaa I_4(u)&=& \iint_{x^{\alpha}y(\frac{\log x}{\log y}+1)\leq u}1\,dxdy
\leq \iint_{x^{\alpha}y\frac{\log x}{\log y}\leq u}1\,dxdy
=\iint_{\frac{y}{\log y}\leq \frac{u}{x^{\alpha}{\log x}}}1\,dxdy \\[4pt]
&\leq&C\int_1^{(u/\log u)^{1/\alpha}}
\int_{y\leq\frac{u\log u}{x^{\alpha}{\log x}}
-\frac{\alpha u}{x^{\alpha}}-\frac{u \log\log x}{x^\alpha \log x}}
 1\,dydx\\[4pt]
&=&C\int^{(u/\log u)^{1/\alpha}}
\Big(\frac{u\log u}{x^{\alpha}{\log x}}
-\frac{\alpha u}{x^{\alpha}}-\frac{u \log\log x}{x^\alpha \log x}\Big)dx
\\[4pt]
&=& C\bigg(\frac{u\log u}{1-\alpha}
\bigg\{\Big[\frac{x^{1-\alpha}}{\log x}\Big]^{(u/\log u)^{1/\alpha}}
+\int^{(u/\log u)^{1/\alpha}}\frac1{x^{\alpha}(\log x)^2}\,dx\bigg\}\\
&&\hskip1pc-\frac{\alpha u}{1-\alpha}
\Big(\frac{u}{\log u}\Big)^{(1/\alpha)(1-\alpha)}
-\frac{u}{1-\alpha} \Big(\frac{u}{\log u}\Big)^{1-1/\alpha}
\frac{\alpha\log\log u}{\log u}\Big(1+{\cal O}\Big(\frac{\log\log u}{\log u}
\Big)\Big)\bigg)\\
&\le&C\left(\frac{u\log u}{1-\alpha}
\bigg\{\frac{\big(u/\log u\big)^{(1/\alpha)(1-\alpha)}}
{\frac1{\alpha}(\log u-\log\log u)}
+\tilde{C}\,\Big[\frac1{1-\alpha}\cdot
\frac{x^{1-\alpha}}{(\log x)^2}\Big]^{(u/\log u)^{1/\alpha}}\bigg\}\right.
\\&&\hskip4pc\left.
-\frac{\alpha u}{1-\alpha}\Big(\frac{u}{\log u}\Big)^{(1/\alpha)(1-\alpha)}
-\frac{u}{1-\alpha} \left(\frac{u}{\log u}\right)^{1-1/\alpha}
\frac{\alpha\log\log u}{\log u}\right)\\
&\le&Cu\log u\frac{\big(u/\log u\big)^{(1/\alpha)(1-\alpha)}}{(\log u)^2}
=C\Big(\frac{u}{\log u}\Big)^{1/\alpha}\,.
\eeaa
As for the lower bound,
\beaa
I_4(u)&=&\iint_{x^{\alpha}y\frac{\log x+\log y}{\log y}\leq u}1\,dxdy
\geq\iint_{x^{\alpha}y\frac{\log u+\log u}{\log y}\leq u}1\,dxdy
=\int^{u/\log u}\int_{x\leq \Big(\frac{u\log y}{2y\log u}\Big)^{1/\alpha}}
1\,dxdy\\[4pt]
&=&\Big(\frac{u}{\log u}\Big)^{1/\alpha}\int^{u/\log u}
\Big(\frac{\log y}{2y}\Big)^{1/\alpha}\,dy
\ge C\Big(\frac{u}{\log u}\Big)^{1/\alpha}\,,
\eeaa
since the last integral has a limit as $u\to \infty$.
Hence, $I_4(u)\approx (u/\log u)^{1/\alpha}$ and
$M(x)=(x/\log x)^{1/\alpha}$.

\subsubsection*{Proof of (\ref{UnglMG})} For the general case the
corresponding integral can be estimated from above by \beaa
I_5(u)&=& \int\int_{\frac{xy (r(x)+r(y))}{L_1(x)L_2(y)}\le u} 1\,
dx\, dy
\leq \int\int_{\frac{y}{L_2(y)}\le u\, \frac{L_1(x)}{x\,r(x)}}1\,dy\,dx\\[5pt]
&\leq & C\, u\, L_2(u) \int^{u L_1(u)/r(u)}  \frac{L_1(x)}{x\,r(x)}\,dx
\le C\, u\,\frac{L_2(u)}{r(u)} \int^{u L_1(u)/r(u)} \frac{L_1(x)}{x}\,dx\\
&\le&  C\,u\,\frac{L_2(u)}{r(u)} \int^{u\,L_1(u)/r(u)} \frac{L_1(x)}{x}\,dx\,,
\eeaa
where we used \cite{gjs}, Lemma 3.1, for the second to last
inequality.

Interchanging the roles of the two variables will in general lead to
a poorer estimate.

In the cases discussed in detail the estimates are precise enough,
since the lower bound
\beaa
I_5(u)&=& \int\int_{\frac{xy(r(x)+r(y))}{L_1(x)L_2(y)}\le u}1\,dxdy
\geq\int_{u^\beta}^{c\,u L_1(u)/r(u)}
\int_{\frac{xy(r(u)+r(u))}{L_1(u^\beta)L_2(y)}
\leq u}1\, dx dy\\
&=&\int_{u^\beta}^{c\,u L_1(u)/r(u)} \int_{y\leq \frac{uL_2(u^\beta)}{2r(u)}
\frac{L_1(x)}{x}}\,1\, dy dx \ge
C \frac{uL_2(u^\beta)}{r(u)}\int_{u^\beta}^{u L_1(u)/r(u)}\frac{L_1(x)}{x}\,dx
\eeaa
e.g.\ for $\beta=1/2$, leads to the same functions (up to constants),
that is, we rediscover our results in the examples discussed above.

Finally, since the last integral on the right-hand side increases
faster than the logarithm and since $r(u)=o(\log u)$, it follows that
$M(u)/u \to \infty\,.$

\medskip\noindent {\small Allan Gut, Department of Mathematics,
Uppsala University, Box 480, SE-751\,06 Uppsala, Sweden;\\
Email:\quad \texttt{allan.gut@math.uu.se}\\
URL:\quad \texttt{http://www.math.uu.se/\~{}allan}}
\\[4pt]
{\small Ulrich Stadtm\"uller, Ulm University, Department of Number
Theory and Probability Theory,\\D-89069 Ulm, Germany;\\
Email\quad \texttt{ulrich.stadtmueller@uni-ulm.de}\\
URL:\quad
\texttt{http://www.mathematik.uni-ulm.de/matheIII/members/stadtmueller/stadtmueller.html}}
\end{document}